\documentclass[12pt]{amsart}

\newcommand {\supplus}{\mathop{{\supset}\llap{\raise 
0.5pt\hbox{\normalfont\small+}\hskip 0.5pt}}} 

\newcommand {\subplus}{\mathop{{\subset}\llap{\raise 
0.5pt\hbox{\normalfont\small+}\hskip 0.5pt}}}  

\newcommand {\Kee}    {{\mathbb  K}}

\newcommand {\Ree}    {{\mathbb  R}}

\newcommand {\Zee}    {{\mathbb  Z}}

\newcommand {\fder}   {{\mathfrak{der}}}   %

\newcommand {\fe}     {{\mathfrak{e}}}

\newcommand {\fg}     {{\mathfrak{g}}}    %

\newcommand {\fo}     {{\mathfrak{o}}}

\newcommand {\fvect}  {{\mathfrak{vect}}}   %

\newcommand {\cal} {\mathcal}

\def \opname#1#2%
  {\expandafter\newcommand \csname #1\endcsname {{\mathop{#2}\nolimits}}}

\newcommand{\rmname}[1]
  {\expandafter\newcommand \csname #1\endcsname {{\operatorname{#1}}}}

\newcommand{\rmnameii}[2]
  {\expandafter\newcommand \csname #1\endcsname {{\operatorname{#2}}}}

\rmname{act}
\rmname{Ad}
\rmname{Add}
\rmname{ad}
\rmname{Alt}
\rmname{alt}
\rmname{Ann}
\rmname{antidiag}
\rmname{Ber}
\rmname{ber}
\rmname{Br}
\rmname{card}
\rmname{ch}
\rmname{Char}
\rmname{cem}
\rmname{cj}
\rmname{Cliff}
\rmname{cntr}
\rmname{codim}
\rmname{coind}
\rmname{const}
\rmname{col}
\rmname{cork}
\rmname{cpr}
\rmname{diag}
\rmnameii{Div}{div}
\rmname{Def}
\rmname{Der}
\rmname{Dim}
\rmname{End}
\rmname{Even}
\rmname{Ext}
\rmname{gr}
\rmname{Hom}
\rmname{HT}
\rmnameii{Ht}{ht}
\rmname{hwt}
\rmname{Id}
\rmname{id}
\rmname{ind}
\rmname{Ind}
\rmname{Inf}
\rmname{irr}
\rmname{Le}
\rmname{Lie}
\rmname{lwt} 
\rmname{mult}
\rmname{Mor}
\rmname{nm}
\rmname{Ob}
\rmname{Odd}
\rmname{Osc}
\rmname{per}
\rmname{Pic}
\rmname{pr}
\rmname{pro}
\rmname{Prime}
\rmname{Proj}
\rmname{prt}
\rmname{pt}
\rmname{Q}
\rmname{qet}
\rmname{qtr}
\rmname{rd}
\rmname{rk}
\rmname{row}
\rmname{Res}
\rmname{salt}
\rmname{Sch}
\rmname{SBr}
\rmname{scalar}
\rmname{Ser}
\rmname{sign}
\rmname{Smbl}
\rmname{spin}
\rmname{ssym}
\rmname{str}
\rmname{st}
\rmname{sgn}
\rmname{sq}
\rmname{symm}
\rmname{supp}
\rmname{Supp}
\rmname{St}
\rmname{Spec}
\rmname{Spm}
\rmname{tr}
\rmname{vpt}
\rmname{weyl}
\rmname{Weyl}
\rmname{Witt}

\opname{vvol}  {{v\hspace{-0.1ex}o\hspace{-0.02ex}l\/}}
\opname{pnt}  {\text{\normalfont pt}}
\opname{Span} {{Span}}
\opname{slim} {\overline{\lim}}
\opname{Vol}  {{V\hspace{-0.55ex}o\hspace{-0.02ex}l\/}}
\opname{QVol} {{Q\hspace{-0.3ex}V\hspace{-0.55ex}o\hspace{-0.02ex}l\/}}
\opname{PoVol}{{P\hspace{-0.35ex}o\hspace{-0.25ex}V\hspace{-0.55ex}o\hspace{-0.02ex}l\/}}
\opname{BVol} {{B\hspace{-0.2ex}V\hspace{-0.55ex}o\hspace{-0.02ex}l\/}}
\opname{Par}  {{P\hspace{-0.3ex}a\hspace{-0.05ex}r\/}}

\rmname{Mat}
\rmname{Bil}
\rmname{Diff}
\rmname{Ker}
\rmname{Herm}
\rmname{Coker}
\rmname{Conn}
\rmname{Covect}
\rmname{Vect}
\rmname{Int}

\rmnameii {IM} {Im}
\rmnameii {RE} {Re}

\opname{Aut} {{A\hspace{-0.2ex}u\hspace{-0.1ex}t\/}}
\opname{GL} {{G\hspace{-0.3ex}L}}
\opname{SL} {{S\hspace{-0.3ex}L}}
\opname{Exp} {{E\hspace{-0.2ex}x\hspace{-0.1ex}p\/}}
\opname{GQ} {{G\hspace{-0.2ex}Q}}
\opname{OSp} {{O\hspace{-0.25ex}S\hspace{-0.15ex}p\/}}
\opname{Out} {{O\hspace{-0.25ex}u\hspace{-0.15ex}t\/}}
\opname{Spp} {{S\hspace{-0.2ex}p\/}}
\opname{SpO} {{S\hspace{-0.2ex}p\hspace{-0.02ex}O\/}}
\opname{Pe} {{P\hspace{-0.25ex}e\/}}
\opname{SPe} {{S\hspace{-0.25ex}P\hspace{-0.25ex}e\/}}
\opname{Spin} {{S\hspace{-0.25ex}p\hspace{-0.05ex}i\hspace{-0.1ex}n\/}}
\opname{Iso} {{I\hspace{-0.25ex}s\hspace{-0.1ex}o\/}}
\opname{SSPe} {{S\hspace{-0.25ex}S\hspace{-0.15ex}P\hspace{-0.25ex}e\/}}
\opname{PeU} {{P\hspace{-0.25ex}e\hspace{-0.1ex}U\/}}
\opname{QU} {{Q\hspace{-0.15ex}U\/}}
\opname{U} {{U\/}}

\opname{cGQ} {{\cal G \hspace{-0.2em} Q \/}}
\opname{cSL} {{\cal S \hspace{-0.2em} L \/}}
\opname{cGL} {{\cal G \hspace{-0.2em} L \/}}
\opname{cOSp} {{\cal O \hspace{-0.2em} S \hspace{-0.3em} \it p\/}}
\opname{cPe} {{\cal P \hspace{-1.5pt} \it e\/}}
\opname{cVect} {{\cal V \hspace{-1.5pt} \it e\hspace{-0.1ex}c\hspace{-0.1ex}t\/}}
\opname{cVol} {{\cal V \hspace{-1.5pt} \it o\hspace{-0.1ex}l\/}}
\opname{cAut} {{\cal A \hspace{-0.2em} \it u\hspace{-0.1em}t\/}}
\opname{cCovect} {{\cal C \hspace{-1.5pt}
     \it o\hspace{-0.1ex}v\hspace{-0.1ex}e\hspace{-0.1ex}c\hspace{-0.1ex}t\/}}
\opname{CW} {{C\hspace{-0.15ex}W}}

\opname {Ab}   {{\sf Ab}}
\opname {Alg}   {{\sf Alg}}
\opname {ASch}  {{\sf Aff\;Sch}}
\opname {Funct}   {{\sf Funct}}
\opname {Gr}   {{\sf Gr}}
\opname {Grf}  {{{\sf Gr}_f}}
\opname {Mods}   {{\sf mods}}
\opname {Rings}   {{\sf Rings}}
\opname {Salg}   {{\sf Salg}}
\opname {Sets} {{\sf Sets}}
\opname {SSMan} {{\sf SMan}}
\opname {Top}  {{\sf Top}}
\opname {Vebun}   {{\sf Vebun}}

\newcommand {\tto} {\longrightarrow}

\newcommand {\bcdot}   {\mathbin{\hbox{\raise.4ex\hbox{\bf.}}}} 

\newcommand {\secno} {}
\newcommand {\ssecfont} {\normalfont\bf}

\newtheorem{Theorem}{\secno Theorem}

\newtheorem{Lemma}[Theorem]{\secno Lemma}
\newtheorem{Proposition}[Theorem]{\secno Proposition}

\newenvironment {th*}[1]
    {\gdef\thname{#1} \begin{thn}}%
    {\end{thn}}
\newtheorem{thn}[Theorem] {\thname}

\theoremstyle{definition}

\newenvironment {ex*}[1]
    {\gdef\thname{#1} \begin{exn}}%
    {\end{exn}}
\newtheorem{exn}[Theorem]{\thname}

\theoremstyle{remark}

\newtheorem{Remark}[Theorem]{\secno Remark}

\newenvironment {rem*}[1]
    {\gdef\thname{#1} \begin{remn}}%
    {\end{remn}}
\newtheorem{remn}[Theorem]{\thname}

\newcommand {\ssec}{\subsection*}

\newcommand {\ssbegin}[2]
  {\def \secno {\gdef \secno {}{\ssecfont #1. }}%
   \begin{#2}}
\begin{document}

\title{The Riemann tensor for nonholonomic manifolds}

\author{Dimitry Leites} 

\address{Department of Mathematics, University of Stockholm,
Kr\"aftriket hus 6, SE-106 91, Stockholm, Sweden;
mleites@matematik.su.se}



\begin{abstract} For every {\it nonholonomic} manifold, i.e., manifold
with nonintegrable distribution the analog of the Riemann tensor is
introduced.  It is calculated here for the contact and Engel
structures: for the contact structure it vanishes (another proof of
Darboux's canonical form); for the Engel distribution the target space
of the tensor is of dimension 2.  In particular, the Lie algebra
preserving the Engel distribution is described.

The tensors introduced are interpreted as modifications of the Spencer
cohomology and, as such, provide with a new way to solve partial
differential equations.  Goldschmidt's criteria for formal
integrability (vanishing of certain Spencer cohomology) are only
applicable to ``one half'' of all differential equations, the ones
whose symmetries are induced by point transformations.  Lie's theorem
says that the ``other half'' consists of differential equations whose
symmetries are induced by contact transformations.  Therefore, we can
now extend Goldschmidt's criteria for formal integrability to all
differential equations.
\end{abstract}

\thanks{Partially supported by an NFR grant.  I am thankful to
A.~Vershik and V.~Gershkovich who put the problem; to L.~Faddeev and
A.~Reyman and to R.~Cianci, U.~Bruzzo and E.~Massa for hospitality
when a guest of the Euler IMI in November 1990 and the Genova
University in February 1991, respectively, and where these results
were delivered.  They were main ingredient of the theoretical part of
my ICTP lectures in 1991, cf.  \cite{GL1}.  I thank P.~Grozman who
checked my computations with the help of his SuperLie package
\cite{GL2}.}

\maketitle

\section*{Introduction} 

Let $W_k^n$ be the space of germs of irreducible $k$-dimensional
distributions at the origin of $\Ree^n$.  The orbits of $\Diff
\Ree^n$-action on $W_k^n$ are usually of infinite codimension.  As
shown in \cite{VG1}, there are 3 types of exceptions: (1) $W_k^n$ is a
point, (2) codimension 1 distributions determined as zero set of
$\alpha =0$ for the forms $\alpha$ of maximal class (recall that the
{\it class} of the form is the invariant equal to the least number of
variables actually involved)
$$
\alpha= \begin{cases} dt +\mathop{\sum}\limits_{1\leq i\leq r}
(p_idq_i-q_idp_i) &\text{for $n=2r+1$}\cr 
dt+\mathop{\sum}\limits_{1\leq i\leq r-1} (p_idq_i-q_idp_i) &\text{for
$n=2r$}\end{cases}
$$
called (usually, for  $n$ odd) a {\it contact structure} (see 
\cite{Z}) and (3) the 
{\it Engel distribution} locally given by the system of Pfaff 
equations $\alpha_1 = 0$, $\alpha_2 = 0$, where
$$
\alpha_1 = dx_4 - x_3dx_1,\quad\alpha_2 = dx_3 - x_2dx_1.\eqno{(E)}
$$
Any 4-dimensional manifold with the Engel distribution 
(E) will be called an {\it Engel manifold}.

Both contact and Engel manifolds are examples of {\it nonholonomic}
manifolds (i.e., manifolds with nonintegrable distributions).  Recall
that the term {\it nonholonomic} was coined by H.~Hertz \cite{H} in
his attempts to geometrically describe motions in such a way as to
eliminate the notion of ``force''.  For a historical review see
\cite{VG2}; see also a very interesting paper by Vershik \cite{V} with
first rigorous mathematical formulations of nonholonomic geometry and
indications to various possible applications in various, sometimes
unexpected, areas (like optimal control or macro-economics).  In
\cite{V} Vershik summarizes about 100 years of studies of nonholonomic
geometry (Hertz, Carath\'eodory, Vr\u anceanu, Wagner, Schouten,
Faddeev, Griffiths, Godbillon, to name a few; to his list I'd add that
MathSciNet returns about 700 entries for nonholonomic, and for
synonyms: 70 more for sub-Riemannian,1700 for Finsler, several dozens
for ``cat's problem'' and autoparallel; physics have their own simialr
lists with thousands of entries because there seems to be more,
actually, nonholonomic dynamical systems than holonomic ones; finally
supergravity is also a nonholonomic structure, albeit on {\it
super}manifolds --- ca 1800 more entries).  In \cite{V} Vershik also
declares that ``it is probably impossible to construct an analog of
the Riemann tensor for the general nonholonomic manifold''  though in some
particular cases of small dimension they had been constructed
(Martinet) (and, one can add now, in papers on supergravity).

Here I introduce the invariants that characterize \lq\lq nonflatness"
of nonholonomic manifolds, i.e., positively solve this old problem. 
This definition appeared around 1989 but publication of this paper was
delayed and the definition first (as far as I know) appeared in
contetext of supergravity in \cite{LP}.  For application to SUGRA
(supergravity equations) and other instances where the N-curvature
tensors (the left hand side of SUGRA) are calculated see \cite{GL1},
\cite{LPS}, \cite{LSV}.

The general definition is illustrated with contact and Engel
manifolds.  As expected, for the contact structure these invariants
vanish, thus giving one more proof of the Darboux theorem on the
canonical form of the contact form; one can show that they vanish for
$n$ even as well.  The Engel manifolds are more rigid: there could be
nonflat ones among them.

Observe that our invariants, though natural analogs of the curvature
and torsion tensors, do not coincide on nonholonomic manifolds with
the classical ones and bearing the same name.  Indeed, by definition
on any nonholonomic manifold, there is a nonzero classical torsion
while every contact manifold is flat in {\it our} sense.  To avoid
confusion, we should always add adjective ``nonholonomic'' for the
invariants introduced below.  Since this is too long, we will briefly say
{\it N-torsion}, {\it N-curvature}, etc.

For manifolds with a $G$-structure, the prototypes of the invariants
we consider are well-known: these are {\it structure functions} of the
$G$-structure \cite{St}; an important particular case is $G = O(n)$
when the structure function is called the {\it Riemann tensor}.

A theorem of Lie states that there are two types of differential
equations, considered as submanifolds in the jet bundles: (1) those
whose symmetries are induced by point transformations and (2) those
whose symmetries are induced by contact transformations.  For
equations of the first type Goldschmidt gave criteria for formal
integrability \cite{BCG}.  The notions introduced below ---
generalizations of Spencer cohomology for nonholonomic case ---
provide with a straightforward generalization of Goldschmidt's
criteria for formal integrability to the second half of differential
equations.

Finally, observe that ``non-flatness'' of distributions with
non-linear constraints (such as fields of cones \cite{A}, foreseen by
Vershik in \cite{V}, or fields of spheres, like the ones constructed
by any driver who switched the cruise control ON, can also be tested
by invariants I introduce; in the case of the fields of solids,
however, the Lie algebras involved are of infinite dimension (they
correspond to ``curved Grassmannians'' of submanifolds, rather than
linear subspaces) and feasibility to obtain the final result becomes
very doubtful.

The waves in plasma is another example of a nonholonomic system with
infinite dimensional algebras involved.  This is a pity, since the
curvature tensor is known to be responsible for stability of the
dynamical system expressed by geodesics; so, in principle, we could
just have computed the N-curvature in this case and settle an
intriguing question: ``is cold fusion possible or not?''.

\section*{\S 1.  Structure functions} 

Let $M$ be a manifold of dimension $n$ over a field $\Kee$.  Let
$F(M)$ be the frame bundle over $M$, i.e., the principal $GL(n;
\Kee)$-bundle.  Let $G\subset GL(n; \Kee)$ be a Lie group.  A $G$-{\it
structure} on $M$ is a reduction of the principal $GL(n; \Kee)$-bundle
to the principal $G$-bundle.  Another formulation seems to be clearer:
if one can select transition functions from one coordinate patch to
another so that they belong to $G$ for every intersection pair of
patches, then we have a $G$-structure on $M$.

The simplest $G$-structure is the {\it flat} $G$-structure defined as 
follows.  Let $V$ be $\Kee^n$ with a fixed frame.  Fix a frame in 
$V$.  Having identified $V$ with $T_vV$ via parallel translation along 
$v$, we get a fixed frame in every $T_vV$.  The flat $G$-structure is 
the bundle over $V$ whose fiber over $v\in V$ consists of all frames 
obtained from the fixed one under the $G$-action. 

In textbooks on differential geometry (e.g., in \cite{St}) the obstructions 
to identification of the $k$-th infinitesimal neighborhood of a point 
$m\in M$ on a manifold $M$ with $G$-structure with that of a point of 
the manifold $V$ with the above flat $G$-structure are called {\it 
structure functions} of order $k$.  It is also explained there that 
such an identification is possible provided all structure functions of 
lesser orders vanish.

To formulate the precise statement, we need some definitions.  Let $S^i$ 
denote the operator of the $i$-th symmetric power.  Set
$$
\fg_{-1} = T_mM,\; \;   \fg_0 = \fg = \Lie (G);
$$ 
for any $i > 0$ and $k= 1$, \dots, $i$ set:
$$
\fg_i = \{X\in \Hom(\fg_{-1}, \fg_{i-1})\mid X(v_1)(v_2, v_3, ..., v_i)  = 
X(v_2)(v_1, v_3, ..., v_i)\;\text{ for any }\; v_k\in \fg_{-1}\}.  
$$
Now, set  $(\fg_{-1}, \fg_{0})_* = \mathop{\oplus}\limits_{i\geq -1} \fg_i$. 
	
Suppose that the $\fg_0$-module $\fg_{-1}$ is
faithful. Then, clearly,  
$$(\fg_{-1}, \fg_{0})_*\subset \fvect (n) = \fder~
\Kee[[x_1,\ldots , x_n ]],\; \text{ where}\; n = \dim~ \fg_{-1}.
$$
It is subject to an easy verification that the Lie algebra structure 
on $\fvect (n)$ induces same on $(\fg_{-1}, \fg_{0})_*$.  The Lie 
algebra $(\fg_{-1}, \fg_{0})_*$ will be called the {\it Cartan's 
prolong} (the result of Cartan's {\it prolongation}) of the pair 
$(\fg_{-1}, \fg_{0})$.

Let $\Lambda^i$ be the operator of the $i$-th exterior power; set
$C^{k,s}_{(\fg_{-1}, \fg_0)_*} = \fg_{k-s}\otimes
\Lambda^s(\fg_{-1}^*)$; we usually drop the subscript of $C^{k,
s}_{(\fg_{-1}, \fg_0)_*}$ or indicate only $\fg_0$.  Define the
differential $\partial _s: C^{k, s}\tto C^{k, s+1}$ setting for any
$v_1, \ldots , v_{s+1}\in V$ (as usual, the slot with the hatted
variable is ignored):
$$
( \delta _sf)(v_1, \ldots , v_{s+1}) = \sum _i (-1)^i[f(v_1, \ldots , {\hat
v}_{s+1-i}, \ldots , v_{s+1}), v_{s+1-i}].
$$

As expected, $\delta _s\delta _{s+1} = 0$.  The homology of this 
bicomplex is called {\it Spencer cohomology} of the pair $(\fg_{-1}, 
\fg_0)$ and denoted by $H^{k,s}_{(\fg_{-1}, \fg_0)_*}$.

\begin{Proposition} {\em (\cite{St})} Structure functions constitute the 
space of the $(k,2)$-th Spencer cohomology of the Cartan prolongation 
of the pair $(T_mM, \fg)$ where $\fg = \Lie (G)$.
\end{Proposition}

The Riemann tensor on an $n$-dimensional Riemannian manifold is the
structure function from $H^{k, 2}_{(\fg_{-1}, \fo(n))_*}$.

\begin{Remark} (cf.  \cite{G}.)
$$
\mathop{\oplus}\limits_kH^{k, s}_{(\fg_{-1}, \fg_0)_*}= H^s(\fg_{-1}; (\fg_{-1}, \fg_0)_*).
$$  
\end{Remark}

This remark considerably simplifies calculations, in particular, if
the Lie algebra $(\fg_{-1}, \fg_0)_*$ is simple and finite
dimensional: we can apply the Borel-Weil-Bott- \ldots (BWB) theorem. 
In the nonholonomic case considered in what follows the remark
provides with a compact definition\footnote{Cf.  with the problems
encountered in the pioneer paper \cite{T}, where the case $d=2$ for
Lie algebras was independently introduced.} of structure functions. 
We can recover the bigrading (explicitly preserved in \cite{T}) at any
moment but to work with just one grading is much simpler, as with
generating functions of series instead of individual terms.  Moreover,
there are a few general theorems that help to calculate Lie algebra
cohomology and no theorems on Spencer cohomology.

\section*{\S 2. Structure functions of nonholonomic manifolds}

To embrace contact-like structures we have to slightly generalize the 
notion of Cartan prolongation: the tangent space to a point of a 
nonholonomic manifold is naturally rigged with a nilpotent Lie algebra 
structure (e.g., for any odd dimensional contact manifolds this is the 
Heisenberg algebra), cf. \cite{VG}.

\ssec{2.1.  Nonholonomic manifolds} (\cite{VG1}, \cite{VG2}) Let $M^n$ 
be an $n$-dimensional manifold with a nonintegrable distribution $D$.  
Let
$$
D= D_1\subset D_2 \subset D_3 \dots \subset D_d
$$
be the sequence of strict inclusions, where the fiber of $D_i$ at a 
point $x\in M$ is $$D_{i-1}|_x\cup[D_1, D_{i-1}]|_x
$$
(here $[D_1, D_{i-1}]=\{[X, Y]\mid  X\in D_1, Y\in D_{i-1}\}$) and $d$ is 
the least number such that
$$
D_{d}|_x\cup[D_1,
D_{d}]|_x = D_{d}|_x.
$$
In case $D_d = TM$ the distribution is called {\it completely 
nonholonomic}.  The number $d = d(M)$ is called the {\it 
nonholonomicity degree}.  A manifold $M$ with a distribution $D$ on it 
will be referred to as {\it nonholonomic} one if $d(M)\neq 1$.

Here we will consider nonholonomic manifolds only when they are
completely nonholonomic.  This does not diminish generality, actually,
since if the distribution in question is not completely nonholonomic,
we can take the integral manifold corresponding to it; the restriction
of the same distribution onto this integral manifold is completely
nonholonomic, all right.  If, nevertheless, we {\it have to}, for some
reason, consider the general case it is easy {\it mutatis mutandis};
it corresponds to {\it nontransitive} Shchepochkina prolongs, cf. 
below, and the general machinery is easily adjusted to them).

Let $n_i(x) = \dim D_i|_x$.  The distribution $D$ is called a {\it 
regular} one if all the dimensions $n_i$ are constants on $M$.  In what 
follows we will only consider regular distributions. This is also 
not a restriction if everything is smooth.

The tangent bundle over a nonholonomic manifold $(M, D)$ can be
naturally endowed with a bundle of nilpotent Lie algebras as follows. 
Let $L$ be the Lie algebra of vector fields preserving the
distribution in a neighborhood of some point $x\in M$. Naturally,  
$L$ is filtered and let $\fg$ be the associated graded Lie algebra. Set
$$
\fg_- = \mathop{\oplus}\limits _{-d\leq i\leq -1}\fg_i.
$$
Clearly, $d$, the depth of $\fg_-$, coinsides with the nonholonomicity
degree of $D$, see \cite{VG2}.

An analog of $G$-structure for nonholonomic manifolds is defined to 
be the pair $(\fg_-,\fg_0)$, where $\fg_0$ is a subalgebra of the Lie 
subalgebra which preserves the $\Zee$-grading of $\fg_-$, i.e., $\fg_0 
\subset (\fder~\fg_-)_0$.  If $\fg_0$ is not indicated explicitely we 
assume that $\fg_0 =(\fder~\fg_-)_0$, i.e., is the largest one.

Thus, we see that with every nonholonomic manifold $(M, D)$ a natural
$G$-structure is associated: $\Lie(G)=(\fder~\fg_-)_0$.  But the
structure functions of this $G$-structure do not reflect the
nonholonomic nature of $M$.  Indeed, when we compute the structure
functions (be they ``traditional'' $H^{k, 2}_{(\fg_{-1}, \fg_0)_*}$,
or introduced above $H^s(\fg_{-1}; (\fg_{-1}, \fg_0)_*)$) we use the
fact that $\fg_{-1}$ is a commutative Lie algebra (partial derivatives
commute).

To take the nonholonomic nature of $M$ into account we need an analog
of the above Proposition for the case when the natural basis of the
tangent space consists not of partial derivatives but rather of
covariant derivatives corresponding to the connection determined by
the same Pfaff forms whose vanishing determines the distribution.
Therefore, instead of $\fg_{-1}$ we have a nilpotent algebra
$\fg_{-}$.  

So we have to define 

(1) the simplest nonholonomic structure --- the ``flat'' one, and

(2) the analog of Cartan prolong when $\fg_{-1}$ is replaced with
$\fg_{-}$,

(3) (we have already understood) what is the analog of $\fg_{0}$ if it
is not mentioned as is usually the case (only the distribution is
given), but it remains to figure out

(4) what is the analog of $(\fg_{-1}, \fg_0)_*)$, and, finally, 

(5) what is the analog of $H^{k, 2}_{(\fg_{-1}, \fg_0)_*}$.

Let us answer these remaining questions.

Fix a frame $F$ in $\Kee^n$.  Having identified $T_vV$ with $V$ by
means of the transformation that preserves $D$,  we fix a
frame --- the image of $F$ --- in each $T_vV$.  The nonholonomic
structure of the nonholonomic manifold $(\Kee^n, D)$ is the pair of
bundles: the frame bundle and the distribution $D$.  The nonholonomic
structure of $(\Kee^n, D)$ will be considered to be {\it flat} if the
frames over $v$ are obtained from the fixed one by means of the
$G$-action, where $G$ is a group whose Lie algebra is
$(\fder~\fg_-)_0$.

\ssec{2.2.  Structure functions of a nonholonomic manifold} Given a
pair $(\fg_-,\fg_0)$ as above, define its $k$-th {\it Shchepochkina
prolong} (first published by I.~Shchepochkina \cite{Sh} for Lie
superalgebras; see also more accessible \cite{Sh5}) for $k> 0$ to be:
$$
\fg_k = (i(S^*(\fg_-)^*\otimes \fg_0)\cap j(S^*(\fg_-)^*\otimes \fg_-))_k, 
$$
where the subscript singles out the component of degree $k$ and where 
$i: S^{k+1}(\fg_{-1})^*\otimes \fg_{-1}\tto S^{k}(\fg_{-1})^*\otimes 
\fg_{-1}^*\otimes\fg_{-1}$ is a natural embedding and $j: 
S^{k}(\fg_{-1})^*\otimes \fg_{0}\tto S^{k}(\fg_{-1})^*\otimes 
\fg_{-1}^*\otimes\fg_{-1}$ a natural map.

Similarly to the case when $\fg_-$ is commutative, define $(\fg_-,
\fg_0)_*$ to be $\mathop{\oplus}\limits_{k\geq -d} \fg_k$; then, as is
easy to verify, $(\fg_-, \fg_0)_*$ is a Lie algebra.  By the same
reasons as for the manifolds with $G$-structure (\cite{St}), the
space $H^2(\fg_-; (\fg_-, \fg_0)_*)$ is the space of structure
functions --- obstructions to the identification of the infinitesimal
neighborhood of a point on the manifold with a nonholonomic structure
(given by $\fg_-$ and $\fg_0$) with same on a flat nonholonomic
manifold with the same $\fg_-$ and $\fg_0$.

The space $H^2(\fg_-; (\fg_-, \fg_0)_*)$ naturally splits into
homogeneous components whose degree will be called the {\it order} of
the structure function; the minimal order of structure function is
equal to $2-d$.  As in the case of a commutative $\fg_-=\fg_{-1}$, the
structure functions of order $k$ can be interpreted as obstructions to
flatness of the nonholonomic manifold with the $(\fg_-,
\fg_0)$-structure provided the obstructions of lesser orders vanish. 
Observe that for a nonholonomic manifold the order of structure
functions is no more in direct relation with the numbers of the
infinitesimal neighborhoods of the points we wish to identify:
distinct partial derivatives bear different \lq\lq weights".

When I told this to A.~Beilinson in 1989 he winced and, as I got this
far, he said that though he conceeded the idea is of interest,
he doubted it can work in the general case: all my constructions are
graded, whereas with every distribution only filtered algebras $\fg_-$
are naturally associated.  This incredulity is justified but
fortunately, different filtrations on $\fg_-$ with the same graded
$\fg_-$ are governed precisely by the structure functions of {\it
lesser} orders and all of them must vanish for the N-curvature to be
well-defined. So, given a filtered $\fg_-$, just take the associated 
graded algebra and proceed as indicated above. This indicates 
a way to study the stable normal forms of Pfaff equations listed in 
\cite{Z}.

\section*{\S 3.  The Lie algebra $\fe$ of Engel vector fields} 

Let us describe the Lie algebra $\fe$ of vector fields that preserve 
the Engel distribution.  Since each of the forms $\alpha_1$ and 
$\alpha_2$ is a familiar contact form, we know enough to demonstrate 
(we leave it as an easy exercise for the reader) the following 
statement.

\begin{Lemma} A weight basis of the $\Zee$-graded Lie algebra 
$\fe$ is as follows (indicated are the degrees): 
\noindent
$$
\begin{tabular}{|c|c|c|c|c}
\hline
$-3$&$-2$&$-1$&$0$&$n\geq1$\cr
\hline
$\partial_4$&$\partial_3
+x_1\partial_4$&$\partial_1$&$
x_1\partial_1+x_2\partial_2+2x_3\partial_3+3x_4\partial_4$&\cr
&&$\partial_2+x_1\partial_3+\frac{1}{2}x_1^2\partial_4$&$x_1\partial_1-
x_2\partial_2+x_4\partial_4$&$
x_1^{n+1}\partial_2+\frac{x_1^{n+2}}{n+2}\partial_3+
\frac{x_1^{n+3}}{n+3}\partial_4$\cr
&&&$x_1\partial_2+\frac{x_1^{2}}{2}\partial_3+
\frac{x_1^{3}}{3}\partial_4$&\cr
\hline
\end{tabular}
$$
\end{Lemma}

The structure of the Engel algebra $\fe$ seem to have never been 
described; so let us do it.  First, introduce a shorthand for some 
elements of $\fe$:
$$
\begin{array}{rl}
    D_3&=\partial_3+x_1\partial_4;\\
E&=x_1\partial_1+x_2\partial_2+2x_3\partial_3+3x_4\partial_4;\\
H&=x_1\partial_1-x_2\partial_2+x_4\partial_4;\\
X_n&=x_1^{n+1}\partial_2+\frac{x_1^{n+2}}{n+2}\partial_3+
\frac{x_1^{n+3}}{n+3}\partial_4\; \; \; \text{  for }\ n\geq -1.
\end{array}
$$

Now, it is clear that $\fe$ is solvable with
$$
\fe/[\fe, \fe] = \Span(E, H);\quad 
\fe^{(1)}/[\fe^{(1)}, \fe^{(1)}] = \Span(\partial _1), 
\text{ where }\fe^{(1)}=[\fe, \fe].
$$

The negative part, $\fe_{-}$ plays for $\fe$ the same role the
Heisenberg algebra plays for the contact algebra.  This negative part
was, of course calculated earlier, see, e.g., \cite{G1} where it also
bears the name of Engel algebra; perhaps one should say the ``small''
and the ``total'' ones.  It seems very strange that the total algebra
had never beed desribed; I am sure such a description is buried in
works of classics.

In review to \cite{G1} [MR 96k:58010] it is written ``There is a
Darboux theorem for Engel structures: any two are locally
diffeomorphic''.  Vershik also quoted this statement as a well-known. 
But I can not place it.  Indeed, by Darboux theorem for the contact
structure I understand the statement that says that the contact form
can be locally reduced to a canonical form , cf.  \cite{Z}.  (In
contradistinction, a related statement that the nondegenerate
differential 2-form $\omega$ can only be reduced to a canonical form
if and only if $d\omega=0$, cf.  \cite{P}, so the statement on
reducibility of the nondegenerate {\it closed} 2-form, sometimes
encountered in textbooks on Differential Geometry, is a tautology.)

Concerning Engel structures, there is, for example, paper \cite{G2}
with examples of several non-equivalent Engel structures on $\Ree^4$. 
Doesn't it contradict the statement on an analog of Darboux theorem? 
So let us calculate the N-curvature of the contact and Engel
structures: to find out if these invariants vanish, hence, there are
analogs of Darboux theorem for them.

\section*{\S 4.  Structure functions of the contact and Engel structures} 

\'Eli Cartan worked before ``(co)homology'' was defined, still, he
introduced a notion that helps to calculate some of these
cohomologies.  This notion is {\it involutive Lie algebra of vector
fields}.  Serre gave a nice reformulation of this notion (\cite{St});
thanks to Serre's statement for the case when $(\fg_-, \fg_0)_*$ is a
simple vectorial Lie algebra instead of computing $H^2(\fg_{-1};
(\fg_{-1}, \fg_0)_*)=\mathop{\oplus}\limits_kH^{k, 2}_{(\fg_{-1},
\fg_0)_*}$, it only suffices to compute the summand with the least
$k=1$: only this term might be nonzero because $(\fg_-, \fg_0)_*$ is
involutive.  This really helps, see examples in \cite{LPS}.

Regrettably, Serre's theorem is only formulated so far for algebras of
depth 1.  So I calculated $H^2(\fg_{-}; (\fg_{-}, \fg_0)_*)$ by hands;
later the result was confirmed by independent computer calculations up
to convincingly high degree:

\ssbegin{4.1}{Theorem} The space of nonholonomic structure functions on any
Engel manifold is of dimension $2$.  All structure functions are of
order $2$, their representatives are:
$$
\partial _1^*\wedge
X_{-1}^*\otimes H+ \partial _1^*\wedge
X_{-1}^*\otimes E\ \text{  and}\; \partial _4^*\wedge
D_{3}^*\otimes \partial _4.
$$
The first of these cycles is \lq\lq pure" in the sence that it does
not represent its cohomology class modulo any coboundary, the
coboundaries vanish.  The other one is mixed: one can add any
coboundary to it.
\end{Theorem} 

The following is another form of the classical theorem (see
\cite{LP}):

\ssbegin{4.2}{Theorem} {\em (Darboux theorem)} The space of
nonholonomic structure functions on any contact manifold vanishes.
\end{Theorem}

 \end{document}